\documentclass[a4paper,reqno, 10pt]{amsart}

\usepackage[margin=2.5cm,papersize={19.75cm,28cm}]{geometry}

\usepackage{amssymb}           
\usepackage{hyperref}
\usepackage{eucal}
\usepackage{epsfig}
\usepackage{amsmath} 
\usepackage{color}
\newtheorem{definition}{Definition}[section]

\newtheorem{remarkth}[definition]{Remark}
\newtheorem{exampleth}[definition]{Example}
\newenvironment{remark}{\begin{remarkth}\upshape}{\hfill$\diamond$\end{remarkth}}

\renewcommand{\emph}[1]{{\bfseries\itshape{#1}}}

\makeatletter

\begin{document}

\title[Time-minimum control of quantum purity for 2-level Lindblad equations]{Time-minimum control of quantum purity for 2-level Lindblad equations}
\author[W. Clark, A. Bloch, L. Colombo, P. Rooney]{William Clark, Anthony Bloch, Leonardo Colombo, Patrick Rooney}

\begin{abstract}
We study time-minimum optimal control for a class of quantum two-dimensional dissipative systems whose dynamics are governed by the Lindblad equation and where control inputs acts only in the Hamiltonian. The dynamics of the control system are analyzed as a bi-linear control system on the Bloch ball after a decoupling of such dynamics into intra- and inter-unitary orbits. The (singular) control problem consists of finding a trajectory of the state variables solving a radial equation in the minimum amount of time, starting at the completely mixed state and ending at the state with the maximum achievable purity. 

The boundary value problem determined by the time-minimum singular optimal control problem is studied  numerically. If controls are unbounded, simulations show that multiple local minimal solutions might exist. To find the unique globally minimal solution, we must repeat the algorithm for various initial conditions and find the best solution out of all of the candidates. If controls are bounded, optimal controls are given by  bang-bang controls using the Pontryagin minimum principle. Using a switching map we construct optimal solutions consisting of singular arcs. If controls are bounded, the analysis of our model also implies classical analysis done previously for this problem.
\end{abstract}

\maketitle


\centerline{\it Dedicated to J\"urgen Scheurle}

\section{Introduction}
 Control of quantum conservative (Hamiltonian) systems has been widely studied in the last few decades from both theoretical and interdisciplinary points of view, see for example \cite{BCS}, \cite{2-3level}, \cite{Jover}, \cite{Brocket1}, \cite{Brocket2}. Recently, there has been a growing interest in control of open (dissipative, non-Hamiltonian) quantum systems because of their applications to physics, chemistry and quantum computing. For example, there has been interest in the control of the rotation of a molecule in the gaseous phase by using laser fields in dissipative media \cite{RS}, where the dissipation is due to molecular collisions. 
Other applications include control of  spin dynamics by magnetic fields in nuclear magnetic resonance \cite{EB}, as well as the construction of quantum computers \cite{RB}.

The aim of this work is to study time-minimal optimal control of two-level quantum systems in a  dissipative environment, where we assume that the dissipation is Markovian (the dynamics depend only on the present state and not its history) and time-independent.  In this case the evolution of the density matrix of the system can be described by a quantum dynamical semi-group and the Lindblad master equation \cite{Altafini}, \cite{Breuer}, \cite{lindblad1976}. The boundary value problem determined by the time-minimum singular optimal control problem is studied  numerically with unbounded and bounded controls respectively. A preliminary version of this work, focusing mainly on the energy-minimization problem with unbounded control, appears in the conference paper \cite{CDC}. We extend the results of \cite{CDC} to the time-minimum optimal control problem including bounded controls. 
 
We study time-minimal controls with $\lVert u\rVert_{\infty}\leq 1$.  
In our feedback classification of trajectories, firstly we reproduce the results in \cite{BC2003}, \cite{5288564}, since their system fits in those studied in this work for particular choices of parameters. We also complete the study done in \cite{BC2003} and \cite{5288564}, where controls are constant and unbounded.
 
The state space for a \textit{closed finite-dimensional quantum system} is an $n$-dimensional projective Hilbert space, $P(\mathcal{H})$. Typically one drops the requirement of projectiveness, and instead work with unit vectors. To preserve the length of the vectors, \textit{time evolution is unitary} $\displaystyle{U(t_1,t_2)|\psi(t_1)\rangle = |\psi(t_2)\rangle}$, and the evolution is described by the \textit{Schr\"{o}dinger equation}
\begin{equation}
\frac{d}{dt}|\psi(t)\rangle = -iH(t)|\psi(t)\rangle,
\end{equation}
where $H$ is the Hermitian Hamiltonian. Here the \textit{ket}-bracket describes the vector associated with an observable state.






The \textit{density operator}, $\rho$, describes a probabilistic ensemble of states. It is given by a positive semi-definite Hermitian operator $\rho$ with $\text{Tr}(\rho)=1$ and $\text{Tr}(\rho^2)\leq 1$. The \textit{purity} of a density operator describes how close $\rho$ is to a single state. It is usually defined as
$\displaystyle{P_2(\rho)={\text{Tr}(\rho^2)}\in [1/n,1]}$, where the unique operator that has a purity of $1/n$ is $\frac{1}{n}I_{n\times n}$, called the \textit{completely mixed state.}  The dynamics for the purity operator is described by the \textit{von Neumann equation}
\begin{equation}
\frac{d}{dt}\rho = [-iH,\rho].
\end{equation}
Notice that this dynamical system preserves the purity of $\rho$ (because the system is iso-spectral). A consequence is that if the quantum system is controlled by its Hamiltonian, there is no control over its purity, since one cannot directly alter the eigenvalues.


\textit{Open quantum systems} are quite different. For such systems, dissipation occurs when we allow the system to interact with the environment. The full picture is an integro-differential equation called the Nakajima-Zwanzig (NZ) equation. To make the dissipation purely of differential form, one usually makes two assumptions: the dissipation is Markovian (i.e. the dissipation only depends upon the current state, not past history) and the dissipation is time-invariant.

Under these assumptions, the dynamics of the density operator is given by the \textit{Lindblad master equation} \cite{lindblad1976}
$$\frac{d}{dt}\rho = [-iH,\rho]+\sum_{j=1}^{N}\left(
L_j\rho L_j^{\dagger} - \frac{1}{2}\left\{ L_j^{\dagger}L_j,\rho\right\}\right)
$$ where the $L_j$ are called the Lindblad operators, $N$ denotes the number of Lindblad operators, and $\{\cdot,\cdot\}$ is the anti-commutator: $\{A,B\}=AB+BA$. We assume the Lindblad operators are traceless.

The structure of the paper is as follow: Section $2$ introduces the Lindblad equation and the reformulation of the equations as a first-order dynamical system on the unit (Bloch) ball. The time-minimal control problem studied in this work is formulated after defining the apogee and escape chimney. Section $3$ is devoted to the study of the control dynamics when there is only one Lindblad term. The case of constant controls is also analyzed. Section $4$ studies the time-minimum optimal control problem. The resulting boundary value problem is studied  numerically in two scenarios: unbounded controls and bounded controls.
\section{Statement of the problem}
\subsection{The Lindblad equation}
\label{section2} 
An open quantum system is described by a density operator $\rho$, which is a trace-one positive semi-definite Hermitian operator on an $n$-dimensional complex Hilbert space $\mathcal{H}$. If the dissipation is Markovian and time independent, the density operator obeys the Lindblad equation
(see \cite{Breuer} and \cite{lindblad1976} for details) 
\begin{equation}\label{eq:master}
\frac{d\rho}{dt} = [-iH,\rho]+\sum_{j=1}^{N} L_j\rho L_j^{\dagger} - \frac{1}{2}\{L_j^{\dagger}L_j,\rho\},
\end{equation} where $N$ denotes the number of Lindblad operators,
$[\cdot,\cdot]$ denotes commutator of matrices, $H$ is the Hermitian Hamiltonian, $\dagger$ represents the Hermitian transpose and $\left\{L_j\right\}$ are the Lindblad operators. The purity of the system is defined as $P_{2}(\rho)={\text{Tr}(\rho^2)}$. We assume that controls appear in the Hamiltonian $H$, and not in the Lindblad operator  $\displaystyle{\mathcal{L}_{D}:=\sum_{j=1}^{N} L_j\rho L_j^{\dagger} - \frac{1}{2}\{L_j^{\dagger}L_j,\rho\}}$.



\indent When $n=2$, the density operator can be identified with a vector in the Bloch ball (the unitary sphere with its interior) \cite{PhysRevA.93.063424}, \cite{RoBlRe}. In this special case, we can change  from the operator point-of-view to that of the Bloch ball by considering $\rho=\frac{1}{2}\left(I+q_1\sigma_1+q_2\sigma_2+q_3\sigma_3\right)$ where $q\in S^{2}$, $I$ is the $2\times 2$ identity matrix, and $\sigma_{j}$ are the Pauli matrices.

We can reformulate the Lindblad equation (\ref{eq:master}) as a first-order dynamical system on the unit ball. Using the derivation given in \cite{PhysRevA.93.063424} (see appendix A therein),  equation \eqref{eq:master} is equivalent to
\begin{equation}\label{eq:bloch}
\frac{d\vec{q}}{dt} = \vec{b} + (A-\text{tr}(A))\vec{q} + \vec{u}\times\vec{q}, \hbox{ where } A :=\frac{1}{2}\sum_j \vec{l_j}\overline{\vec{l}_j}^T+\overline{\vec{l}_j}\vec{l}_j^T,\hbox{ and }
\vec{b} := i\sum_j\! \vec{l}_j \times \overline{\vec{l}_j},
\end{equation} where the bar represents the complex conjugate of matrices.
The vectors $\vec{l}_j$, $\vec{u}$ are the traceless parts of $L_j$ and $H$ respectively with $H = \displaystyle{h_0I+\sum_{k=1}^3 \! u_k\sigma_k}$ where $\sigma_k$ are the Pauli matrices. Notice that the matrix $A$ is positive semi-definite. 
We will define $B:=A-\text{tr}(A)$.

The controlled Hamiltonian dynamics cannot achieve purity one (see \cite{tannor}) and in general cannot affect the purity of the state or transfer the states between unitary orbits. To control purity, one must use the dissipative dynamics to move between orbits as in \cite{PhysRevA.93.063424}, \cite{RoBlRe} and \cite{Rooney2}. 
If $q\neq 0$, equation \eqref{eq:bloch} can be written in terms of the radial component (i.e., $r=\lVert q\rVert$). 

In \cite{PhysRevA.93.063424}  it is shown that the purity $P_{2}(\rho)$ is equal to $\displaystyle{\frac{(1+r^2)}{2}}$, where $r=\lVert q\rVert$.   Knowing that $r^2=\langle q,q\rangle$, we find that
\begin{equation*}
\begin{split}
2r\dot{r} &= 2\langle \dot{q},q\rangle
= 2\langle b+u\times q + Bq,q\rangle = 2\langle b,q\rangle +2\langle Bq,q\rangle
= 2r\langle b,\hat{q}\rangle + 2r^2\langle B\hat{q},\hat{q}\rangle
\end{split}
\end{equation*} where $\hat{q}$ is the unit vector associated with $q$, $\hat{q} = q /||q||$. Therefore, 
\begin{equation}\label{eq:redial}
r\frac{dr}{dt} = \langle q, b+Bq\rangle := f(q).
\end{equation} So, we can control the purity via the orientation of the corresponding unit vector.  Hereafter we refer to $\displaystyle{f(q)}$ as \textit{the purity derivative}.



\subsection{The apogee and the escape chimney} \label{section3}
The purity derivative (\ref{eq:redial}) is independent of the controls $u_j$. It is illuminating to examine the regions in the Bloch ball where the purity derivative is positive. To do this, we examine the zeros of $f$. 

Define two sets, $\mathcal{U} = \{q|f(q)\geq 0\}$ and the ellipsoid $\mathcal{M} = \{q|f(q)=0\}$. To find $\mathcal{M}$, we define a new function $f_q(r) := f(rq)$ where $q\in S^2$. Finding the roots of $f_q$ will let us solve for $\mathcal{M}$ in spherical coordinates. Note that
\begin{equation}\label{eq:radialzeros}
f_q(r) = \langle q,Bq\rangle r^2 + \langle q,b\rangle r.
\end{equation} so, the nonzero root is
\begin{equation}\label{eq:nonzeroroot}
g(q) := -\frac{\langle q,b\rangle}{\langle q,Bq\rangle}.
\end{equation}
Notice that $g$ is always defined, since $B$ is negative-definite. Also note that the maximum of $g$ must be bounded by $1$, so the Bloch ball is invariant under (\ref{eq:bloch}).\\
\indent We define $q_{apogee}$ to be \textit{the apogee} of the ellipsoid $\mathcal{M}$. i.e.
\begin{equation}\label{eq:endpoint}
q_{apogee} := \arg\max_{q\in\mathcal{M}} ~ \lVert q \rVert,
\end{equation}
which can be found by maximizing \eqref{eq:nonzeroroot} on $S^2$. Therefore, the apogee of the ellipsoid $\mathcal{M}$ will be the state with the maximal achievable purity. We call the interior of the ellipsoid $\mathcal{U}$ the \textit{escape chimney}. 

 The control problem consists of finding a trajectory of the state variables solving \eqref{eq:redial}, starting at the completely mixed state (i.e., $r=0$) and ending  at the apogee. 
It is important to note, however, that the dynamics (\ref{eq:redial})  has a singularity at the origin since $$\frac{dr}{dt}=\frac{f(q)}{r}.$$ Also note that the apogee cannot be reached in finite time since it is not possible to reach equilibrium points in finite time. Since the purity derivative is independent of controls, the purity obeys an autonomous first-order dynamical system that cannot reach its fixed point. 

By considering the boundary conditions representing initial and final states
\begin{equation}\label{eq:startandend}
q_0 = \varepsilon \vec{b}/\lVert b\rVert, \hspace{0.15in} q_f = (1-\delta)q_{apogee}
\end{equation}
with $\varepsilon>0$ and $\delta>0$ sufficiently small, the optimal control problem can be stated as follows:

\textbf{Problem:} Let $J:S^2\times U\rightarrow\mathbb{R}$ be a cost functional dependent on the state as well as the controls. The time-minimum optimal control problem consists of finding a control, $u^{*}:[0,t_f]\rightarrow U$ satisfying the dynamics (\ref{eq:redial}) such that $q(0)=q_0$, $q(t_f)=q_f$, and $\displaystyle{u^{*} = \arg\min\! \int_0^{t_f} dt.}$


\section{Dynamics for Lindblad equations}
In the special case when $N=1$ in (\ref{eq:master}), $\vec{b}$ becomes an eigenvector of $B$. This fact lets us simultaneously diagonalize $B$ and rotate $\vec{b}$ into the first coordinate. By additionally taking $u_1=u_2=0$, the third component of $q$ in equation \eqref{eq:bloch} becomes uncontrolled and exponentially decays to zero. Dropping this coordinate, our system collapses to a two-dimensional underactuated bi-linear control system:
\begin{equation}\label{eq:twodimenions}
\begin{array}{rcl}
\dot{x} &=& b_1 + \alpha_1 x -uy \\
\dot{y} &=& b_2 + \alpha_2 y + ux
\end{array} 
\end{equation} where $q=(x,y)$, $\vec{b}=(b_1,b_2)$, $u=u_3$, $\alpha_1,\alpha_2 <0 $ are the coefficients of the matrix $A$,  and the maximum of \eqref{eq:radialzeros} is bounded by $1$. We assume without loss of generality that $\alpha_2<\alpha_1$.

\begin{remark}
System \eqref{eq:twodimenions} reduces in special cases to those  studied in \cite{Bonnarcyber} and \cite{Sugny}. Our system \eqref{eq:twodimenions} is of the form $\dot{q}=F(q)+uG(q)$, $F(q)=(b_1+\alpha_1x,b_2+\alpha_2y)^{T}$ and $G(q)=(-y,x)^{T}$ while in \cite{Bonnarcyber} and \cite{Sugny}, after the change of coordinates $x=y$, $y=z$, the drift term and control vector field are $F(q)=(-\gamma x,\gamma_{-}-\gamma_{+}y)$ and $G(q)=(-y,x)$ with parameters $\Gamma, \gamma_{+},\gamma_{-}$. Therefore, by setting $b_1=0, \alpha_1=-\Gamma, b_2=\gamma_{-}, \alpha_2=-\gamma_{+}$ in \eqref{eq:twodimenions} we recover the situation in \cite{Bonnarcyber} and \cite{Sugny}. Note also that the constraints in \cite{Bonnarcyber} and \cite{Sugny} become $-2\alpha_1\geq-\alpha_2\geq|b_2|$, which always holds in our case as long as $\alpha_2<\alpha_1<0$.
\end{remark} 

Note that in the absence of controls, the (asymptotically stable) fixed point for the drift field is $q^{*}=(q_1^{*},q_2^{*}):=-B^{-1}b$. 
To determine when $q^{*}=q_{f}$ (when the apogee is at $q^{*}$), one can use polar coordinates and study the critical points of $g(\theta)$. Denoting $\theta^{*}:=\hbox{atan2}\left(\frac{q_{2}^{*}}{q_1^{*}}\right)$, and using MATLAB's symbolic toolbox: $$\frac{d}{d\theta}g(\theta)\Big{|}_{\theta=\theta^{*}}=\frac{-b_1b_2(\alpha_1-\alpha_2)\sqrt{(b_1/\alpha_1)^2+(b_1/\alpha_2)^2}}{\alpha_2b_1^2+\alpha_1b_2^2}.$$ Hence, $\theta^{*}$ is a critical point if $\alpha_1=\alpha_2$, or $b_1, b_2=0$. Looking at the sign of the second-derivative for $\alpha_1=\alpha_2$: $$\frac{d^2}{d\theta^2}g(\theta^{*})=\frac{-\sqrt{b_1^2+b_2^2}}{|\alpha_1|}<0,$$ and therefore $g$ achieves its local maximum at $\theta^{*}$.

If $b_1=0$, $$\frac{d^2}{d\theta^2}g(\theta^{*})=\frac{-b_2^2|\alpha_1\alpha_2|(2\alpha_1-\alpha_2)}{\alpha_2^3|\alpha_1b_2|}$$ and therefore $g$ achieves its local maximum at $\theta^{*}$ if $2\alpha_1-\alpha_2<0$. If $b_2=0$, we similarly obtain the same conclusion under the condition $2\alpha_2-\alpha_1<0$. Therefore we conclude that when controls are turned off, we can achieve maximum purity if (1) $\alpha_1=\alpha_2$, (2) $b_1=0, 2\alpha_1-\alpha_2 < 0$, or (3) $b_2=0,2\alpha_2-\alpha_1 <0$. 

While using a \textit{constant control} is not the optimal method, it offers the most simplistic choice of control for our problem. In the absence of controls, the trajectory hits the ellipsoid $\mathcal{M}$ either  to the left or to the right of $q_f$. In the presence of controls the system introduces a \textit{swirl} to the dynamics. We would like to create a constant swirl  such that $q_f$ the end-state of the system.

 Assume we have constant controls $u$. Fixed points for the dynamics determined by  \eqref{eq:twodimenions} are $$q^* = \frac{1}{u^2+\alpha_1\alpha_2}\left\{\begin{array}{c}
-\alpha_2b_1-b_2u\\
-\alpha_1b_2 + b_1u
\end{array}\right\}.$$
Letting $\theta^* = \text{atan2}(q^*_2,q^*_1)$, we want to find $u$ such that $g'(\theta^*)=0$. According to MATLAB's symbolic toolbox, this is equivalent to finding a root for the following cubic polynomial (which is guaranteed to have a real root):
\begin{equation}\label{eq:polynomial}
\lVert b\rVert^2 u^3 + 3b_1b_2(\alpha_2-\alpha_1)u^2 + (2\alpha_2^2b_1^2+2\alpha_1^2b_2^2 - \alpha_1\alpha_2(b_2-b_1))u +\alpha_1\alpha_2b_1b_2(\alpha_1-\alpha_2) = 0.
\end{equation}
Observe that $u=0$ is a solution to (\ref{eq:polynomial}) exactly when $b_1=0$, $b_2=0$, or $\alpha_1=\alpha_2$, which is consistent with what we expect. 

We are interested in studying time-minimum control problems with bounded controls, so a natural question to ask is which conditions are necessary for the root of (\ref{eq:polynomial}) to be bounded by 1?

We rewrite \eqref{eq:polynomial} as $p(u)=au^3+bu^2+cu+d$ (note that $a>0$). We need to first determine how many real, distinct roots there are. Let $\Delta:= 4b^2-12ac$ be the discriminant of $p'(u)$.

\begin{enumerate}
	\item $\Delta \leq0$: Equation \eqref{eq:polynomial} can only have one real root. So necessary conditions for the root to be bounded by 1 is for $\text{sgn}(p(-1)/p(1))=-1$, or $p(1)=0$. In this case, we require:
	$$\frac{-a+b-c+d}{a+b+c+d} \leq 0$$ or 
	$$a+b+c+d  = 0.$$
	\item $\Delta > 0$: Equation \eqref{eq:polynomial} can have 1, 2, or 3 real roots. If the above condition is satisfied, then we have an odd number of roots in $[-1,1]$, but it will not detect an even number of roots. Assume there are two roots, $u_1, u_2 \in [-1,1]$. Then, by the mean value theorem, there must exist $4u^*\in (u_1,u_2)$ such that $p'(u^*)=0$ (which can be found by the quadratic formula). 
\end{enumerate}

So the conditions for a root in $[-1,1]$ are
	$$a+b+c+d = 0$$
	or
	$$\frac{-a+b-c+d}{a+b+c+d} \leq 0$$
	or
	$$\frac{-a+b-c+d}{a+b+c+d} > 0,~~\text{and}~~\min |u_{\pm}| < 1$$ where 
	$$u_{\pm} = \frac{-2b\pm \sqrt{\Delta}}{6a}.$$

\section{Time-minimal optimal control problem}
In this section we study the time-minimum optimal control problem. The resulting boundary value problem is studied numerically. We study two different scenarios: unbounded and bounded controls. If controls are unbounded, simulations show that multiple local minimal solutions might exist. If controls are bounded, optimal controls are given by  bang-bang controls which are obtained by studying the associated adjoint system to the Pontryagin minimum principle.

\subsection{Time-minimal optimal control problem with unbounded controls}
We want to find (unbounded) controls that steer (\ref{eq:twodimenions}) with end points (\ref{eq:startandend}) in the minimal amount of time possible. That is, we want to find a minimal solution to the functional
\begin{equation}\label{eq:2dlagrangian}
\min \int_{0}^{t_f} \! dt = \min \int_{x_0}^{x_f} \! \frac{dt}{dr}\frac{dr}{dx} \, dx, 
\end{equation} where $x(0)=x_0$ and $x(t_f)=x_f$.
The derivative $dr/dt$ is found in $(\ref{eq:redial})$, and $dr/dx = x+yy'$. So we wish to minimize a functional with integrand 
\begin{equation}\label{eq:lagrangian}
I(q,q') = \int\! L(q,q')  dt = \int_{x_0}^{x_f} \! \frac{x+yy'}{\langle q,b+Bq\rangle} \, dx.
\end{equation}
This Lagrangian is not hyperregular, so the Euler-Lagrange equations will fail to yield meaningful results \cite{baillieul2008nonholonomic}.
To get around this problem, we implement the Rayleigh-Ritz numerical algorithm \cite{hoffman2001numerical}. We assume $y(x)$ is a sum of linearly independent functions
\begin{equation}\label{eq:yform}
y(x) = y^0(x) + \sum_{i=1}^M c_iy^i(x),
\end{equation}
where $y^0(x_0) = y_0$, $y^0(x_f)=y_f$, and $y^i(x_0)=y^i(x_f)=0$. Specifically, we will take the following functions as a basis of polynomials for our approximation:
\begin{equation}\label{eq:guessfunctions}
y^i(x) = (x-x_0)(x-x_f)^i,\quad y^0(x) = \frac{y_f-y_0}{x_f-x_0}(x-x_0)+y_0,\quad i=1,\ldots,M
\end{equation} with $M$ an arbitrary integer.
To minimize the functional \eqref{eq:lagrangian}, the following $M$ equations must hold: 
\begin{equation}\label{eq:partials}
I_{c_i} = \int_{x_0}^{x_f}\! \frac{\partial}{\partial c_i} L \, dx = 0.
\end{equation}  where $\displaystyle{I_{c_i}=\frac{\partial I}{\partial c_i}}.$
This can be done by symbolically computing $\partial L/ \partial c_i$ and numerically integrating, using, for instance, a $4^{th}$ order Runge-Kutta method. In order to find the optimal values to the $c_i$'s, we construct the function $\nu:\mathbb{R}^M\to\mathbb{R}$ with:
$$
\nu(c)=\left( \sum_{i=1}^M\left(I_{c_i}\right)^2 \right)^{1/2}.$$
\subsubsection{Numerical results:}
We will work an example with parameter values $b_1=1$, $b_2=2$, $a_1=-3$ and $a_2=-4$. Additionally, we will take $\varepsilon = \delta = 10^{-3}$. \\
\indent Solving for the apogee \eqref{eq:endpoint} in polar coordinates, we get $q_{apogee} = [0.4079,0.4493]^T$. Figure $1$ shows a simulation of the trajectory for the $7^{th}$ order curve.

\begin{table}[ht]
\centering 
\begin{tabular}{c c c}
\hline                       
M & Time & Energy\\ [0.5ex]
\hline                  
1 & 1.9372  &  7.5830  \\
3 & 1.9366 & 8.6873  \\
5 & 1.9361 & 1.6368  \\
7 & 1.9359 & 1.3765 \\ [1ex]      
\hline
\end{tabular}\label{table:2d}
\caption{Numerical results from time-minimal controls with solutions of various orders.}

\end{table}

\begin{figure}[h!]
	\includegraphics[scale=0.35]{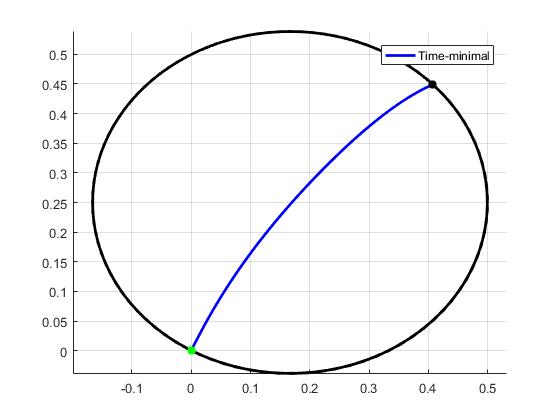} \includegraphics[scale=0.35]{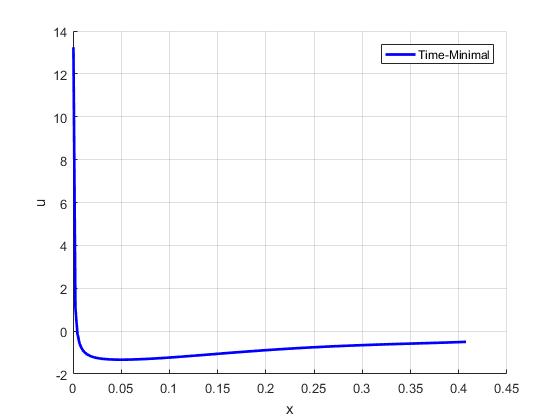}

	\caption{Trajectory and controls of the $7^{th}$ order curve. The black ellipse is the escape chimney.}
\end{figure}

Since all the computations are done independent of time, we report a plot of $u$ versus $x$ in Figure $1$.


An expected  downside of the simulations is that multiple local minimal solutions might exist. To find the unique global minimal solution, we repeat the algorithm for various initial conditions. We then find the best solution out of all of the candidates. For this case, we ran the algorithm 25 times and the initial conditions were uniformly randomly chosen in the $l^{\infty}$-ball of radius 2.

\subsection{Time-minimal optimal control problem with bounded controls}
The goal now is to perform time-minimal controls with $\lVert u\rVert_{\infty}\leq 1$. In our feedback classification of trajectories, we firstly reproduce the results in \cite{BC2003}, \cite{5288564}, as a special case of our model by performing the time-minimal local syntheses \cite{Boscain}. We additionally present an example of our model that is not included in \cite{BC2003} or \cite{5288564}.


To find the optimal control, we use the Pontryagin minimum principle (\cite{kirk2004optimal} Section $5.4$). Consider the two-dimensional underactuated bi-linear control system
\begin{equation}\label{eq:apprev}
\dot{q}=b+Bq+g(q)u,\hspace{0.15in}g(q)=\left\{ \begin{array}{c}
-y\\ x \end{array}\right\}
\end{equation}
From the Pontryagin minimum principle, the Hamiltonian is
\begin{equation}\label{eq:Hamiltonian}
\mathcal{H}(q,p,u) = 1+\langle p, b+Bq+g(q)u\rangle.
\end{equation}
Here we denote $q=(x,y)$, $f(q)=b+Bq=(b_1+\alpha_1x,b_2+\alpha_2y)^{T}$ and $g(q)=(-y,x)^{T}$. Following the analysis performed in \cite{kirk2004optimal}, the optimal controls are given by bang-bang controls as
\begin{equation}\label{eq:controls}
u = \begin{cases}
1 & \langle p, g(q) \rangle <0 \\
-1 & \langle p,g(q) \rangle >0\\
\text{Undefined} & \langle p,g(q) \rangle = 0.
\end{cases}
\end{equation}
To represent trajectories with switching controls we follow the method of time-minimal local syntheses used in \cite{BC2003} and \cite{Boscain}, which avoids using the adjoint system $\displaystyle{\dot{p}=-p\left(B+u\frac{\partial g}{\partial q}\right)}$, $u=\pm 1$, provided by Pontryagin principle. The main idea of the construction provided in \cite{BC2003} (Section $4.5$) is to construct the optimal control as the concatenation of arcs $\sigma_{+}$, $\sigma_{-}$ corresponding to $u=1$, $u=-1$, and singular arcs. The switching function is given by $\Phi(t)=\langle p(t), g(q(t))\rangle$ and the switching takes place when $\Phi(t)=0$. Note that $g$ rotates $\pi$ degrees when we follow an arc $\sigma_{+}$ or $\sigma_{-}$.

Assume $t_0=0$, $t_1=t>0$ are two consecutive switching times where the controls are $u=\varepsilon$, $\varepsilon=\pm 1$. Then by definition of $\Phi$, we must have
$$\langle p(0),g(q(0))\rangle = \langle p(t),g(q(t))\rangle=0.$$
As in \cite{BC2003} (Section $4.5.5$), let $v(\cdot)$ be the solution of the variational equation $\displaystyle{\dot{v}=\left(B+u\frac{\partial g}{\partial q}\right)v}$, $u=\pm1$, such that $v(t)=g(q(t))$, where we integrate backwards from time $t_1$ to $t_0$. By construction, $\langle p(0),v(0)\rangle = 0$. Thus $p(0)$ is orthogonal to $g(q(0))$ and $v(0)$, and therefore $v(0)$ and $g(q(0))$ are collinear. In other words, we want
$$\det (g(q(0)),v(0))=0.$$
As it was shown in \cite{Boscain} (Chapter $2$), defining $\theta(t)$ to be the angle between $g(q(0))$ and $v(0)$, the switching occurs when $\theta(t)=0$ modulo $\pi$ due to reflections along the vertical axis and one has
\begin{equation}\label{eq:solvevar}
v(0)=e^{-t\,\text{ad}(f+\varepsilon g)}g(q(t)).
\end{equation}
In the analytic case, this can be solved by a power series. 

To explicitly compute \eqref{eq:solvevar}, we lift our system to the semi-direct product Lie group $G=GL_2(\mathbb{R})\ltimes \mathbb{R}^2$. Its corresponding semi-direct product Lie algebra, $\mathfrak{g} = \mathfrak{gl}_2(\mathbb{R})\ltimes \mathbb{R}^2$, has a Lie bracket defined by
$$[(A,a),(B,b)] = ([A,B],Ab-Ba)$$ 
for $A,B\in  \mathfrak{gl}_2(\mathbb{R})$ and $a,b\in \mathbb{R}^2$.

We want to compute $\text{ad}(F+\varepsilon G)$ where $F=(B,b)$ and $G=(C,0)$ with $C(q)=g(q)$. 
To compute $\text{ad}(F+\varepsilon G)$, we take the canonical basis for $\mathfrak{gl}_2(\mathbb{R})\ltimes \mathbb{R}^2$ denoted by $\{e_i\}_{i=1}^{6}$ where $$e_1=\left( \begin{array} {cc}
1 & 0 \\ 0 & 0 \end{array}\right), \quad e_2=\left( \begin{array} {cc}
 0 & 1 \\ 0 & 0 \end{array}\right),\quad e_3=\left( \begin{array} {cc}
0 & 0 \\ 1 & 0 \end{array}\right),\quad e_4=\left( \begin{array} {cc}
0 & 0 \\ 0 & 1 \end{array}\right),$$
$e_5=(1,0)^{T}$ and $e_6=(0,1)^{T}$.

 Let $Z = F+\varepsilon G$. Computing the brackets, we get:
\begin{equation}\label{eq:brackets}
\begin{split}
[Z,e_1] &=  (\varepsilon e_2 + \varepsilon e_3 ,-b_1 e_5) \\
[Z,e_2] &= (-\varepsilon e_1 +(\alpha_1-\alpha_2)e_2 +\varepsilon e_4, -b_2 e_5)\\
[Z,e_3] &= (-\varepsilon e_1 +(\alpha_2-\alpha_1)e_3 + \varepsilon e_4 ,-b_1 e_6)\\
[Z,e_4] &= (-\varepsilon e_2 -\varepsilon e_3 -,b_2 e_6) \\
[Z,e_5] &= (0,\alpha_1 e_5 +\varepsilon e_6) \\
[Z,e_6] &= (0,-\varepsilon e_5 + \alpha_2 e_6). \\
\end{split}
\end{equation}
Compiling \eqref{eq:brackets}, we get the adjoint representation
\begin{equation}\label{eq:bigmatrix}
\text{ad}(Z) = \left( \begin{array}{cccccc}
0 & \varepsilon & \varepsilon & 0 & -b_1 & 0 \\
-\varepsilon & \alpha_1-\alpha_2 & 0 & \varepsilon & -b_2 & 0 \\
-\varepsilon & 0 & \alpha_2-\alpha_1 & \varepsilon & 0 & -b_1 \\
0 & -\varepsilon & -\varepsilon & 0 & 0 & -b_2 \\
0 & 0 & 0 & 0 & \alpha_1 & \varepsilon \\
0 & 0 & 0 & 0 & -\varepsilon & \alpha_2
\end{array}\right)^T.
\end{equation}


We want to find $e^{-t\,\text{ad}(Z)}(-e_2+e_3)$. Because this involes exponentiating a $6$ by $6$ matrix, we used MATLAB. To find the collinearity and the switching times, we want to solve for $\text{det}(G(q_0),e^{-t\text{ad}(Z)}(-e_2+e_3)(q_0))=0$.
Using numerical methods, we find $t$ such that 
$$\det (G(q_0),e^{-t\,\text{ad}(Z)}(-e_2+e_3)(q_0)) = 0$$
and switch $\varepsilon$ to $-\varepsilon$ at that point.

\begin{figure}[h!]
	\includegraphics[scale=0.2]{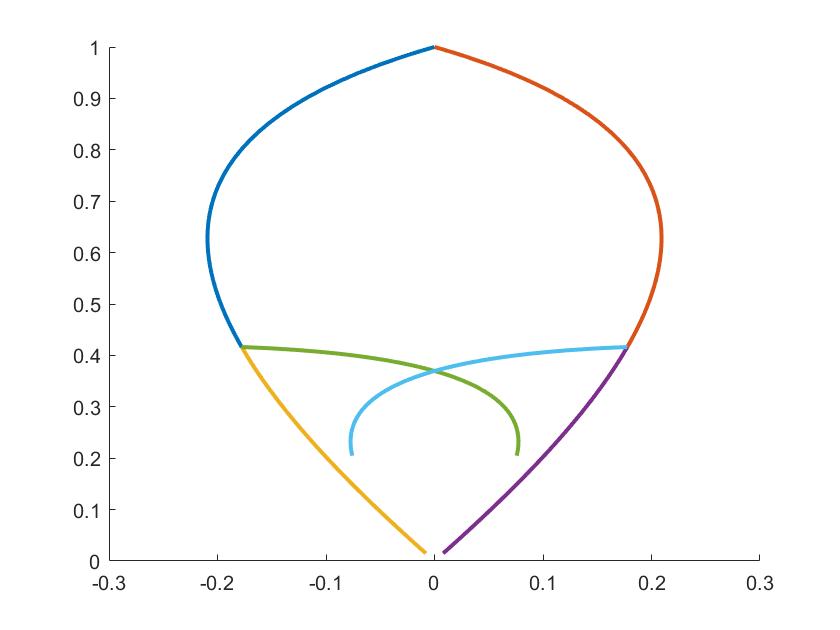} \includegraphics[scale=0.2]{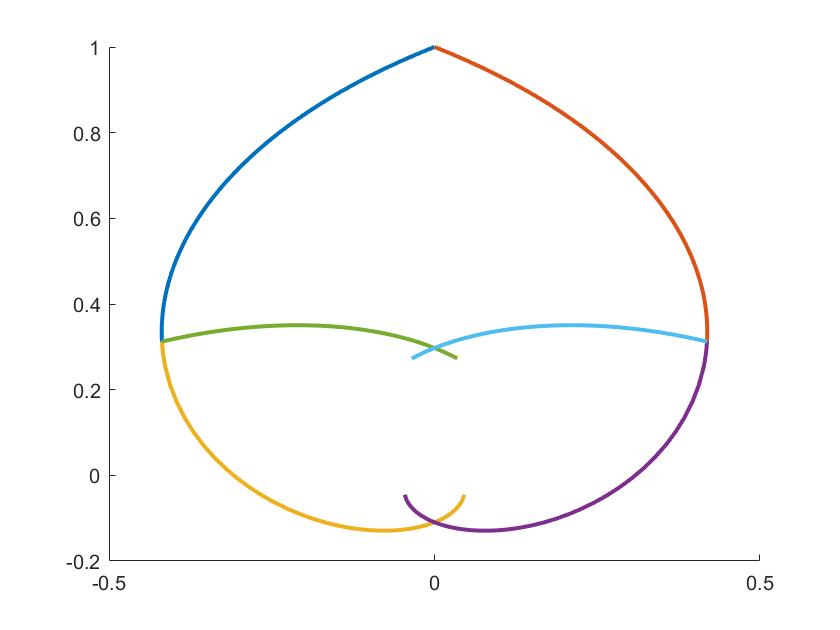}
	\caption{Left: Example where $b=[0;0]$, $\alpha=-3$, and $\beta=-0.6$. Right: Example where $b=[0;0]$, $\alpha=-0.8$, and $\beta=-0.6$.}
	\label{figure2}
\end{figure}

In Figure \ref{figure2} we show a simulation for a specific choice of parameters to mimic the results given in \cite{BC2003} and \cite{5288564}. In the figure, the blue curves represent a control $u=1$ and the red curve if $u=-1$. The green curve and light blue curve are switching signs of the controls (as predicted by the above analysis). Finally, the yellow curve and purple curve corresponds to the case when one never switched signs and one continues with the original controls.

\begin{figure}[h!]
	\includegraphics[scale=0.4]{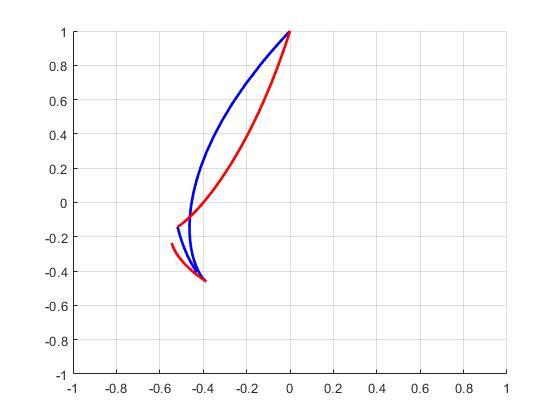}
	\caption{Example where $b=[-2;-1]$, $\alpha = -4$, and $\beta = -3$. The red curves show the trajectory when $u=-1$ and blue when $u=1$.}
	\label{figure3}
\end{figure}
In the general case where $b\ne 0$ (see Figure \ref{figure3}), the switching times are no longer symmetric. The switching times for this example are
\begin{table}[h!]
	\centering
	\begin{tabular}{c|c|c}
		Initial control & Time until first switch & Time between first and second switches \\ \hline
		$u = 1$ & 2.1943 & 0.4685 \\
		$u = -1$ & 1.1532 & 0.4905 \\
	\end{tabular}
	\caption{Switching times for the example in Figure \ref{figure3}.}
	\label{table:table1}
\end{table}
\section{Final discussion}
We studied time-minimum global optimal control problems for dissipative open quantum systems where the dynamics is described by the Lindblad equation and controls are both unbounded and bounded. We have transformed such a control system into a bi-linear singular control system in the Bloch ball and have come up with the construction of a numerical algorithm to design optimal paths to achieve a desired point given initial states close to the origin of the Bloch ball in both optimal control problems. 

We would like to extend our results to higher-order dimensional systems, as well as study energy-minimum global optimal control problems for dissipative open quantum systems where the controls are bounded. It would also be interesting to determine the best basis of functions for the Rayleigh-Ritz methods as well as the best order of solutions to use.

{\bf Acknowledgments:} This research was supported in part by grants 
from the National Science Foundation and Air Force Office of Scientific
Research.


\addtolength{\textheight}{-12cm}   







\begin{thebibliography}{1}

\bibitem{Altafini}
C. Altafini
\newblock {\em Controllability properties for finite-dimensional quantum Markovian master equations}.
\newblock {\em J. Math. Phys.} vol.44, p. 2357, 2002.

\bibitem{baillieul2008nonholonomic}
A.~Bloch, J.~Baillieul, P.~Crouch, and J.~Marsden.
\newblock {\em Nonholonomic Mechanics and Control}.
\newblock Interdisciplinary Applied Mathematics. Springer New York, 2008.

\bibitem{BC2003}
B.~Bonnard and M.~Chyba.
\newblock Singular trajectories and their role in control theory. 
\newblock {\em Math\'ematiques $\&$ Appl.}, \textbf{40} Springer-Verlag, Berlin,  2003.

\bibitem{5288564}
B.~Bonnard, M.~Chyba, and D.~Sugny.
\newblock Time-minimal control of dissipative two-level quantum systems: The
  generic case.
\newblock {\em IEEE Transactions on Automatic Control}, 54(11):2598--2610, Nov
  2009.






\bibitem{Bonnarcyber} B. Bonnard and D. Sugny.
\newblock Time-minimal control of dissipative two-level quantum systems: The integrable case
\newblock {\em Control and cybernetics},vol. 38, no. 4A, pp. 1053-1080, 2009.



\bibitem{Boscain} 
U. Boscain, B. Piccoli. \textit{Optimal Syntheses for Control Systems on 2-D Manifolds}, Springer SMAI, VOL.43. (2004)

\bibitem{BCS}
U. Boscain, M. Caponigro, M. Sigalotti.
\newblock Multi-input Schr\"odinger equation: controllability, tracking, and application to the quantum angular momentum.
\newblock {\em J. of Diff. Eqns.}, \textbf{256} (5), 3524-3551, 2014.

\bibitem{Breuer} H. Breuer and F. Pertuccione.
\newblock The theory of open quantum systems. 
\newblock Oxford University Press, 2007.

\bibitem{Brockett1975}
R.~W. Brockett.
\newblock {\em On the Reachable Set for Bilinear Systems}, pages 54--63.
\newblock Springer Berlin Heidelberg, Berlin, Heidelberg, 1975.

\bibitem{2-3level}
G. Charlot, J.P Gauthier,  U. Boscain,  S. Gu\'erin and H. Jauslin.
\newblock Optimal control in laser-induced population transfer for two- and three-level quantum systems.
\newblock {\em J. Math. Phys.}, \textbf{43} (5), 2107-2132, 2002.


\bibitem{CDC}
W. Clark, A. Bloch, L. Colombo, P. Rooney
\newblock Optimal Control of Quantum Purity for $n=2$ Systems. IEEE CDC Conference Proceedings 2017, Melbourne, Australia. To appear.
\bibitem{EB}
R. Ernst, G. Bodenhausen, A. Wokaun.
\newblock {\em Principles of Nuclear Magnetic Resonance in One and Two dimensions}.
\newblock Clarendon, Oxford, 1987.

\bibitem{hoffman2001numerical}
J.D. Hoffman and S.~Frankel.
\newblock {\em Numerical Methods for Engineers and Scientists, Second
  Edition,}.
\newblock{Taylor \& Francis, 2001.}

\bibitem{Jover} J. Jover Galtier. \textit{Open quantum systems: geometric description,
dynamics and control}. Ph.D thesis, Universidad de Zaragoza, Spain. 2017. 

\bibitem{Brocket1}
N. Khaneja, R. Brockett and S. Glaser.
\newblock {\em Time optimal control of spin systems}.
\newblock  {\em Phys. Rev. A}, \textbf{63} 032308, 2001.

\bibitem{Brocket2}
N. Khaneja, S. Glaser and R. Brockett.
\newblock {\em Sub-Riemannian geometry and time optimal control of three spin systems: Coherence transfer and quantum gates}.
\newblock  {\em Phys. Rev. A}, \textbf{65} 032301, 2002.

\bibitem{kirk2004optimal}
D.E. Kirk.
\newblock {\em Optimal Control Theory: An Introduction}.
\newblock Dover Books on Electrical Engineering Series. Dover Publications,
  2004.

\bibitem{lindblad1976}
G.~Lindblad.
\newblock On the generators of quantum dynamical semigroups.
\newblock {\em Comm. Math. Phys.}, 48(2):119--130, 1976.

\bibitem{RS}
S. Ramakrishna, T. Seideman.
\newblock Intense laser alignment in dissipative media as a route of solvent dynamics.
\newblock {\em Phys. Rev Lett.} vol 95, p.113001, 2005. 

\bibitem{RB}
C.~Rangan and P. Bucksbaum.
\newblock Optimality shaped terahertz pulses for phase retrieval in a Rydberg-atom data registrer.
\newblock {\em Phys. Rev. A},89 (18): 188301, 2002.

\bibitem{PhysRevA.93.063424}
P.~Rooney, A.~Bloch, and C.~Rangan.
\newblock Flag-based control of quantum purity for $n=2$ systems.
\newblock {\em Phys. Rev. A}, 93:063424, 2016.

\bibitem{RoBlRe}
P.~Rooney, A.~Bloch, and C.~Rangan.
\newblock Projector-based Control of Orbit Dynamics in
Quantum Lindblad Systems (2017, to appear in the 
IEEE Trans. Aut. Control).
\newblock {\em arXiv:1201.0399v1}, Preprint 2012.
\bibitem{Rooney2}
P.~Rooney, A.~Bloch, and C.~Rangan.
\newblock Decoherence control and purification of two-dimensional quantum density matrices under Lindblad dissipation.
\newblock {\em arXiv:1201.0399v1}, Preprint 2012.

\bibitem{Sugny} D. Sugny, C. Kontz and H.R. Jauslin.
\newblock Time-optimal control of two-level dissipative quantum system.
\newblock Phys. Rev. A, vol. 76, 2007, 023419.

\bibitem{tannor}
D. Tannor and A. Bartana.
\newblock On the interplay of control fields and spontaneous emission in laser cooling.
\newblock {\em J. Phys Chem A} 103: 10359-10363, 1999.

\end{thebibliography}
\end{document}